\documentclass{amsart}
\usepackage{amssymb,latexsym,enumerate,amscd}
\newtheorem{Thm}[equation]{Theorem}
\newtheorem{Prop}[equation]{Proposition}
\newtheorem{Lem}[equation]{Lemma}
\newtheorem{Cor}[equation]{Corollary}

\theoremstyle{remark}

\theoremstyle{definition}

\newtheorem{Def}[equation]{Definition}
\newtheorem*{Not*}{Notation}
\newtheorem*{Def*}{Definition}
\numberwithin{equation}{section}

\DeclareMathOperator{\tr}{tr}

\DeclareMathOperator{\lcm}{lcm}

\DeclareMathOperator{\Mat}{M}
\DeclareMathOperator{\GL}{GL}
\DeclareMathOperator{\diag}{diag}
\DeclareMathOperator{\rationals}{\mathbb{Q}}
\DeclareMathOperator{\complexs}{\mathbb{C}}

\begin{document}

\date{%
Tue Nov 20 12:41:28 EST 2007}

\title[The Atiyah Conjecture]
{The Atiyah conjecture and Artinian rings}

\author[Peter A. Linnell]{Peter A. Linnell}
\address{Department of Mathematics \\
Virginia Tech \\
Blacksburg \\
VA 24061-0123 \\
USA}
\email{linnell@math.vt.edu}
\urladdr{http://www.math.vt.edu/people/plinnell/}

\author[Thomas Schick]{Thomas Schick}
\address{Mathematisches Institut \\
Georg-August-Universit\"at G\"ottingen\\
Bunsenstr.~3 \\
D-37073 G\"ottingen \\
Germany}
\email{schick@uni-math.gwdg.de}
\urladdr{http://www.uni-math.gwdg.de/schick/}

\begin{abstract}
Let $G$ be a group such that its finite subgroups have bounded
order, let $d$ denote the lowest common multiple
of the orders of the finite subgroups
of $G$, and let $K$ be a subfield of $\mathbb{C}$ that is
closed under complex conjugation.  Let $\mathcal{U}(G)$
denote the algebra of unbounded operators affiliated to
the group von Neumann algebra $\mathcal{N}(G)$, and let
$\mathcal{D}(KG,\mathcal{U}(G))$
denote the division closure of $KG$ in $\mathcal{U}(G)$; thus
$\mathcal{D}(KG,\mathcal{U}(G))$ is the smallest subring of
$\mathcal{U}(G)$ containing $KG$ that is closed under taking
inverses.  Suppose $n$ is a positive integer,
and $\alpha \in \Mat_n(KG)$.  Then $\alpha$ induces a bounded
linear map $\alpha \colon \ell^2(G)^n \to \ell^2(G)^n$,
and $\ker\alpha$
has a well-defined von Neumann dimension $\dim_{\mathcal{N}(G)}
(\ker\alpha)$.  This is a nonnegative real number, and one version
of the Atiyah conjecture states that
$d \dim_{\mathcal{N}(G)}(\ker\alpha) \in
\mathbb{Z}$.  Assuming this conjecture, we shall prove that if
$G$ has no nontrivial finite normal subgroup, then
$\mathcal{D}(KG,\mathcal{U}(G))$ is a $d \times d$ matrix ring
over a skew field.  We shall also consider the case
when $G$ has a nontrivial finite normal subgroup, and other subrings
of $\mathcal{U}(G)$ that contain $KG$.
\end{abstract}

\keywords{Atiyah conjecture, group von Neumann algebra}

\subjclass[2000]{Primary: 16S34; Secondary: 20C07, 22D25, 46L99}

\maketitle

\section{Introduction} \label{Sintroduction}

In this paper $\mathbb{N}$ will denote the positive integers
$\{1,2,\dots \}$, all rings will have a 1, subrings will have the
same 1, and if $n \in \mathbb{N}$, then $\Mat_n(R)$ will indicate
the $n \times n$ matrices over the ring $R$ and $\GL_n(R)$ the
invertible matrices in $\Mat_n(R)$.
Let $G$ be a group, let $\ell^2(G)$ denote the Hilbert space with
orthonormal basis the elements of $G$, and let $\mathcal{B}(\ell^2(G))$ denote
the bounded linear operators on $\ell^2(G)$.  Thus we can write
elements $a \in \ell^2(G)$ in the form $\sum_{g\in G} a_gg$, where
$a_g \in \mathbb{C}$ and $\sum_{g\in G} |a_g|^2 < \infty$.
Then $\mathbb{C}G$ acts faithfully on the left of
$\ell^2(G)$ as bounded linear operators via the left regular representation,
so we may consider
$\mathbb{C}G$ as a subalgebra of $\mathcal{B}(\ell^2(G))$.  The weak
closure of $\mathbb{C}G$ in $\mathcal{B}(\ell^2(G))$ is the group
von Neumann algebra $\mathcal{N}(G)$ of $G$.  Also if $n \in
\mathbb{N}$, then $\Mat_n(\mathbb{C}G)$ acts as bounded linear
operators on $\ell^2(G)^n$ and the weak closure of this ring in
$\mathcal{B}(\ell^2(G)^n)$ is $\Mat_n(\mathcal{N}(G))$.
Let 1 indicate the
element of $\ell^2(G)$ which is 1 at the identity of $G$ and zero
elsewhere.  Then the map $\theta \mapsto \theta 1 \colon
\mathcal{N}(G) \to \ell^2(G)$ is an injection, so we may regard
$\mathcal{N}(G)$ as a subspace of $\ell^2(G)$.  We can now define
$\tr \colon \mathcal{N}(G) \to \mathbb{C}$ by $\tr (a) =
a_1$.  For $\alpha \in \Mat_n(\mathcal{N}(G))$, we can extend
this definition by setting
$\tr(\alpha) = \sum_{i=1}^n \tr(\alpha_{ii})$, where $\alpha_{ij}$ are the
entries of $\alpha$.  A useful property is that if $\alpha$
is a positive operator, then $\tr(\alpha) \ge 0$.  Also we can use tr
to give any right
$\mathcal{N}(G)$-module $M$ a well defined dimension
$\dim_{\mathcal{N}(G)} M$, which in general is a non-negative real
number or $\infty$ \cite[\S 6.1]{Lueck02}.  If $e$ is a projection in
$\Mat_n(\mathcal{N}(G))$, then $\dim_{\mathcal{N}(G)}
e\Mat_n(\mathcal{N}(G)) = \tr(e)$.  Furthermore if $\alpha \in
\Mat_n(\mathcal{N}(G))$, so $\alpha$ is a Hilbert space map
$\ell^2(G)^n \to \ell^2(G)^n$, then since $\ell^2(G)^n$ is a right
$\mathcal{N}(G)$-module, $\dim_{\mathcal{N}(G)} \ker\alpha$ is well
defined and is equal to $\dim_{\mathcal{N}(G)} \{\beta \in
\mathcal{N}(G) \mid \alpha\beta = 0\}$.  Finally $\mathcal{N}(G)$ has
an involution which sends an operator to its adjoint; if $a = \sum_{g
\in G} a_g g$, then $a^* = \sum_{g \in G} \overline{a}_g g^{-1}$,
where the bar indicates complex conjugation.

A ring $R$ is called regular, or sometimes von Neumann regular, if
for every $x \in R$, there exists an idempotent $e \in R$ with $xR
= eR$ \cite[Theorem 1.1]{Goodearl91}.  It is called finite, or
directly finite, if $xy = 1$ implies $yx=1$ for all $x,y \in R$.
Finally a $*$-regular ring $R$ is a
regular ring with an involution $^*$ with the property that $x \in R$
and $x^*x = 0$ implies $x = 0$.  In a $*$-regular ring, given $x\in
R$, there is a unique projection $e$ such that $xR = eR$; so $e = e^*
= e^2$.

Let $\mathcal{U}(G)$ denote
the algebra of unbounded operators on $\ell^2(G)$ affiliated to
$\mathcal{N}(G)$ \cite[\S 8]{Lueck02}.  Then the involution on
$\mathcal{N}(G)$ extends to an involution on $\mathcal{U}(G)$, and
$\mathcal{U}(G)$ is a
finite $*$-regular algebra.  Also if $M$ is a right
$\mathcal{N}(G)$-module, then $\dim_{\mathcal{N}(G)} M =
\dim_{\mathcal{N}(G)} M \otimes_{\mathcal{N}(G)} \mathcal{U}(G)$; in
particular $\dim_{\mathcal{N}(G)} e\mathcal{U}(G) = \tr(e)$.

For any subring $R$ of the ring $S$,
we let $\mathcal{D}(R,S)$ denote the division closure of $R$ in $S$;
that is the smallest subring of $S$ containing $R$ that is closed
under taking inverses.  In the case $G$ is a group and $K$ is a
subfield of $\mathbb{C}$, we shall set $\mathcal{D}(KG) =
\mathcal{D}(KG,\mathcal{U}(G))$.
For any group $G$, let $\lcm(G)$ indicate the least common multiple
of the orders of the finite subgroups of $G$, and adopt the
convention that $\lcm(G) = \infty$ if the orders of the finite
subgroups of $G$ are unbounded.  One version of the strong
Atiyah conjecture states that if $G$ is a group with $\lcm(G) <
\infty$, then the $L^2$-Betti numbers of every closed manifold with
fundamental group $G$ lie in the abelian group $\frac{1}{\lcm(G)}
\mathbb{Z}$.  This is equivalent to the conjecture that if $n\in
\mathbb{N}$, $A \in \Mat_n(\rationals G)$ and
$\alpha \colon \ell^2(G)^n \to \ell^2(G)^n$
is the map induced by left multiplication by $A$,
then $\lcm(G) \dim_{\mathcal{N}(G)} \ker\alpha \in \mathbb{Z}$
\cite[Lemma 2.2]{Lueck02a}.
In this paper, we shall consider more generally the case when the
coefficient ring is a subfield of $\mathbb{C}$.
\begin{Def}
Let $G$ be a group with $\lcm(G) < \infty$, and let $K$ be a subfield
of $\mathbb{C}$.  We say
that the \emph{strong Atiyah conjecture} holds for $G$ over $K$ if
\begin{equation*}
\lcm(G) \dim_{\mathcal{N}(G)} \ker\alpha \in \mathbb{Z}\qquad\text{for all }
\alpha \in \Mat_n(KG).
\end{equation*}
\end{Def}
This is equivalent to the conjecture that
if $M$ is a finitely presented $KG$-module, then
$\lcm(G)\dim_{\mathcal{N}(G)} M \otimes_{KG} \mathcal{N}(G)
\in \mathbb{Z}$ \cite[Lemma 10.7]{Lueck02}.
Obviously if $G$ satisfies the strong Atiyah conjecture over
$\mathbb{C}$, then $G$ satisfies the strong Atiyah conjecture over
$K$ for all subfields $K$ of $\mathbb{C}$.
The strong Atiyah conjecture over $\mathbb{C}$ is known for large
classes of groups; for example \cite[Theorem 1.5]{Linnell93} tells
us that it is true if $G$ has a normal free subgroup $F$ such that
$G/F$ is an elementary amenable group. If $K$ is the
algebraic closure of $\rationals$ in $\complexs$, it is known for even
larger classes of groups, for example
\cite[Theorem 1.4]{Dodziuk-Linnell-Mathai-Schick-Yates(2001)}
for groups which are
residually torsion-free elementary amenable. The following result is well
known; see for example \cite[Lemma 3]{Schick00}. \nocite{MR1894160}
\begin{Prop} \label{Pfield}
Let $G$ be a torsion-free group (i.e.~$\lcm(G) = 1$) and let $K$ be a
subfield of $\mathbb{C}$.  Then $G$ satisfies the strong Atiyah
conjecture over $K$ if and only if $\mathcal{D}(KG)$
is a skew field.
\end{Prop}
The purpose of this paper is to generalize Proposition \ref{Pfield}.
We will denote the finite conjugate subgroup of the group $G$ by
$\Delta(G)$, and the torsion subgroup of $\Delta(G)$ by
$\Delta^+(G)$ (this is a subgroup, compare \cite[Lemma 19.3]{MR0314951}).  We
shall prove
\begin{Thm} \label{Tmain}
Let $G$ be a group with $d := \lcm(G) < \infty$
and $\Delta^+(G) = 1$, and
let $K$ be a subfield of $\mathbb{C}$ that is closed under complex
conjugation.  Then $G$ satisfies the strong Atiyah conjecture over
$K$ if and only if $\mathcal{D}(KG)$ is a $d\times d$
matrix ring over a skew field.
\end{Thm}

It seems plausible that if $K$ is a subfield of $\mathbb{C}$
which is closed under complex conjugation and $G$ is a group with
$\lcm(G) < \infty$ which satisfies the Atiyah conjecture over $K$,
then $\mathcal{D}(KG)$ is a semisimple Artinian ring.
However we cannot prove this, though
we are able to prove a slightly weaker result, and to state
this we require the following definition.

\begin{Def} \label{De}
Let $R$ be a subring of the ring $S$.  The extended division closure,
$\mathcal{E}(R,S)$, of $R$ in $S$ is the smallest subring of $S$
containing $R$ with the properties
\begin{enumerate}[\normalfont(a)]
\item \label{De1}
If $x\in \mathcal{E}(R,S)$ and $x^{-1} \in S$, then $x \in
\mathcal{E}(R,S)$.

\item \label{De2}
If $x \in \mathcal{E}(R,S)$ and $xS = eS$ where $e$ is a central
idempotent of $S$, then $e \in \mathcal{E}(R,S)$.
\end{enumerate}
\end{Def}
Obviously $\mathcal{E}(R,S) \supseteq \mathcal{D}(R,S)$.
Note that if $\{R_i\}$ is a collection of subrings of $S$ satisfying
\ref{De}\eqref{De1} and \ref{De}\eqref{De2} above, then
$\bigcap_i R_i$ is also a subring of $S$ satisfying
\ref{De}\eqref{De1} and \ref{De}\eqref{De2}, consequently
$\mathcal{E}(R,S)$ is a well defined subring of $S$ containing $R$.
Also if $G$ is a group and $K$ is a subfield of $\mathbb{C}$,
then we write $\mathcal{E}(KG)$ for $\mathcal{E}(KG,\mathcal{U}(G))$.
Observe that, if $G$ is torsion-free and if the strong Atiyah
conjecture holds for $G$
over $K$, then $\mathcal{D}(KG)$ is a division ring, hence
$x\mathcal{D}(KG)=\mathcal{D}(KG)$ for every $0\ne x\in
\mathcal{D}(KG)$ and consequently $\mathcal{E}(KG)=\mathcal{D}(KG)$
in this case. We are tempted to conjecture that this is always the
case. We hope to show in a later paper that this should
follow from a suitable version of the Atiyah conjecture.

We shall prove
\begin{Thm} \label{Textended}
Let $G$ be a group with $\lcm(G) < \infty$, and let $K$ be a subfield
of $\mathbb{C}$ that is closed under complex conjugation.  Suppose
that $G$ satisfies the strong Atiyah conjecture over $K$.  Then
$\mathcal{E}(KG)$ is a semisimple Artinian ring.
\end{Thm}
Thus in particular if $K$ is a subfield of $\mathbb{C}$ that is
closed under complex conjugation and $G$ is a group with
$\lcm(G) < \infty$ which satisfies the strong Atiyah conjecture over
$K$, then $KG$ can be embedded in a semisimple Artinian ring.
Theorem \ref{Textended} follows immediately from the more general
Theorem \ref{Pextended} in Section \ref{Sproofs}.

In Section \ref{sec:embeddings} we will show, somewhat unrelated to the
rest of the paper, that $KG$ can be embedded in a least subring of
$\mathcal{U}(G)$ that is $*$-regular.

\section{Proofs} \label{Sproofs}

Let $R$ be a subring of the ring $S$ and let $C = \{e \in S \mid
e$ is a central idempotent of $S$ and $eS = rS$ for some $r \in
R\}$.  Then we define
\[
\mathcal{C}(R,S) = \sum_{e\in C} eR,
\]
a subring of $S$.  In the case $S = \mathcal{U}(G)$, we write
$\mathcal{C}(R)$ for $\mathcal{C}(R,\mathcal{U}(G))$.
For each ordinal
$\alpha$, define $\mathcal{E}_{\alpha}(R,S)$ as follows:
\begin{itemize}
\item
$\mathcal{E}_0(R,S) = R$;
\item
$\mathcal{E}_{\alpha+1}(R,S) =
\mathcal{D}(\mathcal{C}(\mathcal{E}_{\alpha}(R,S),S),S)$;
\item
$\mathcal{E}_{\alpha}(R,S) = \bigcup_{\beta < \alpha}
\mathcal{E}_{\beta}(R,S)$ if $\alpha$ is a limit ordinal.
\end{itemize}
Then $\mathcal{E}(R,S) = \bigcup_{\alpha} \mathcal{E}_{\alpha}(R,S)$.
Also in the case $R = KG$ where $G$ is a group and $K$ is a subfield
of $\mathbb{C}$, we shall write $\mathcal{E}_{\alpha}(KG)$ for
$\mathcal{E}_{\alpha}(KG,\mathcal{U}(G))$.  If $A \subseteq
\mathbb{R}$, then $\langle A \rangle$ will indicate the additive
subgroup of $\mathbb{R}$ generated by $A$.

\begin{Lem} \label{Le}
Let $G$ be a group, let $R$ be a subring of $\mathcal{U}(G)$, let
$n \in \mathbb{N}$, and let $x \in R$.
Suppose that $x \mathcal{U}(G) = e\mathcal{U}(G)$ where
$e$ is a central idempotent of $\mathcal{U}(G)$.
Then $\langle \dim_{\mathcal{N}(G)} \beta\mathcal{U}(G)^n \mid
\beta \in \Mat_n(R)\rangle =
\langle \dim_{\mathcal{N}(G)} \alpha\mathcal{U}(G)^n \mid
\alpha \in \Mat_n(R+eR)\rangle$.
\end{Lem}
\begin{proof}
Set $E = eI_n$, the diagonal matrix in $\Mat_n(R +
eR)$ that has $e$'s on the main diagonal and zeros elsewhere.
Then $E$ is a central idempotent in $\Mat_n(\mathcal{U}(G))$.
Obviously
\[
\langle \dim_{\mathcal{N}(G)} \beta\mathcal{U}(G)^n \mid
\beta \in \Mat_n(R)\rangle \subseteq
\langle \dim_{\mathcal{N}(G)} \alpha\mathcal{U}(G)^n \mid
\alpha \in \Mat_n(R+eR)\rangle,
\]
so we need to prove the reverse inclusion.  Let $\alpha \in
\Mat_n(R+eR)$ and write $\alpha = \beta + E\gamma$ where
$\beta,\gamma \in \Mat_n(R)$.  Then we have
\[
\dim_{\mathcal{N}(G)} \alpha \mathcal{U}(G)^n =
\dim_{\mathcal{N}(G)} (\beta+\gamma) E \mathcal{U}(G)^n +
\dim_{\mathcal{N}(G)} \beta (1-E)\mathcal{U}(G)^n.
\]
Since $\dim_{\mathcal{N}(G)} \beta (1-E)\mathcal{U}(G)^n =
\dim_{\mathcal{N}(G)} \beta \mathcal{U}(G)^n - \dim_{\mathcal{N}(G)}
\beta E \mathcal{U}(G)^n$, it suffices to prove that
\[
\dim_{\mathcal{N}(G)} E\beta \mathcal{U}(G)^n \in
\langle\dim_{\mathcal{N}(G)} \delta\mathcal{U}(G)^n \mid
\delta \in \Mat_n(R)\rangle
\]
for all $\beta \in \Mat_n(R)$.  But $E\beta\mathcal{U}(G)^n =
\beta (xI_n) \mathcal{U}(G)^n$ and the result follows.
\end{proof}

Lemma \ref{Le} immediately gives the following corollary.

\begin{Cor} \label{Ce}
Let $G$ be a group, let $R$ be a subring of $\mathcal{U}(G)$, and
let $n\in \mathbb{N}$.
Then $\langle \dim_{\mathcal{N}(G)} \alpha\mathcal{U}(G) \mid
\alpha \in \Mat_n(R)\rangle =
\langle \dim_{\mathcal{N}(G)} \alpha\mathcal{U}(G) \mid
\alpha \in \Mat_n(\mathcal{C}(R))\rangle$.
\end{Cor}
\begin{proof}
Let $e_1,\dots,e_m$ be central idempotents of $\mathcal{U}(G)$ such
that for each $i$, there exists $\alpha_i \in R$ with $e_i
\mathcal{U}(G) = \alpha_i \mathcal{U}(G)$.  Then by induction on
$m$, Lemma \ref{Le} tells us that the result is true if $\alpha
\in \Mat_n(R+ e_1R + \dots + e_mR)$.  Since $\Mat_n(\mathcal{C}(R))$
is the union of $\Mat_n(R+ e_1R + \dots + e_mR)$,
the result is proven.
\end{proof}

\begin{Lem} \label{LCramer}
Let $R$ be a subring of the ring $S$, let $n\in \mathbb{N}$, and
let $A \in \Mat_n(\mathcal{D}(R,S))$.  Then there
exist $0 \le m \in \mathbb{Z}$ and $X,Y \in \GL_{m+n}(S)$ such that $X
\diag(A,I_m) Y \in \Mat_{m+n}(R)$.
\end{Lem}
\begin{proof}
This follows from \cite[Proposition 7.1.3 and Exercise 7.1.4]{Cohn85}
and \cite[Proposition 3.4]{Linnell06}.
\end{proof}

\begin{Lem} \label{Lsame}
Let $G$ be a group and let $K$ be a subfield of $\mathbb{C}$.
Then $\langle \dim_{\mathcal{N}(G)}
x\mathcal{U}(G) \mid x \in \Mat_n(KG),\ n\in \mathbb{N}\rangle =
\langle\dim_{\mathcal{N}(G)} x\mathcal{U}(G) \mid x \in
\Mat_n(\mathcal{E}(KG)),\ n\in \mathbb{N}\rangle$.
\end{Lem}
\begin{proof}
Obviously
\[
\hspace*{-5pt}
\langle\dim_{\mathcal{N}(G)} x\mathcal{U}(G) \mid x \in \Mat_n(KG),
\ n\in \mathbb{N}\rangle \subseteq \langle\dim_{\mathcal{N}(G)}
x\mathcal{U}(G) \mid x \in
\Mat_n(\mathcal{E}(KG)), \ n\in \mathbb{N}\rangle.
\]
We shall prove the reverse inclusion by
transfinite induction.  So let $n\in \mathbb{N}$ and
$x \in \Mat_n(\mathcal{E}(KG))$.
Then we may choose the least ordinal $\alpha$
such that $x \in \Mat_n(\mathcal{E}_{\alpha}(KG))$.  Clearly $\alpha$
is not a limit ordinal, and the result is true if $\alpha = 0$,
so we may write $\alpha = \beta + 1$ for some ordinal $\beta$ and
assume that the result is true for all
$y \in \Mat_n(\mathcal{E}_{\beta}
(KG))$.  By Corollary \ref{Ce} the result is true for all $y \in
\Mat_n(\mathcal{C}(\mathcal{E}_{\beta}(KG)))$ and now the result
follows from Lemma \ref{LCramer}.
\end{proof}

The following result from \cite{Linnell07} will be crucial
for our work here.  Because of this, and because we use a slightly
different formulation, we state it here.

\begin{Lem} \label{Luseful} \cite[Lemma 2]{Linnell07}
Let $G$ be a group, let $n\in \mathbb{N}$, and let $\alpha_1,
\dots,\alpha_n \in \mathcal{U}(G)$.  Then
$(\sum_{j=1}^n\alpha_j\alpha_j^*)\mathcal{U}(G) \supseteq
\alpha_1\mathcal{U}(G)$.
\end{Lem}
\begin{proof}
By induction on $n$ and \cite[Lemma 2]{Linnell07}, we see that
$(\sum_{j=1}^n\alpha_j\alpha_j^*)\mathcal{U}(G) \supseteq
\alpha_1\alpha_1^*\mathcal{U}(G)$.  The result now follows by
applying \cite[Lemma 2]{Linnell07} in the case $\beta = 0$.
\end{proof}

\begin{Thm} \label{Pextended}
Let $G$ be a group and let $K$ be a subfield of $\mathbb{C}$ which is
closed under complex conjugation.
Suppose there is an $\ell \in \mathbb{N}$ such that
$\ell \dim_{\mathcal{N}(G)} \alpha \mathcal{U}(G)^n\in \mathbb{Z}$ for
all $\alpha \in \Mat_n(KG)$ and for all $n \in \mathbb{N}$.  Then
$\mathcal{E}(KG)$ is a semisimple Artinian ring.
\end{Thm}
\begin{proof}
First observe that Lemma \ref{Lsame} tells us that
\begin{equation} \label{Esamecor}
\ell \dim_{\mathcal{N}(G)} \alpha \mathcal{U}(G) \in \mathbb{Z}
\text{ for all } \alpha \in \mathcal{E}(KG).
\end{equation}
Next note that the hypothesis tells us that $\mathcal{E}(KG)$ has
at most $\ell$ primitive central idempotents.  Indeed if $e_1,\dots,
e_{\ell +1}$ are (nonzero distinct) primitive central idempotents,
then $e_ie_j = 0$ for $i\ne j$ and we see that the sum
$\bigoplus_{i=1}^{\ell+1} e_i\mathcal{U}(G)$ is direct.  But
\[
\dim_{\mathcal{N}(G)} \bigoplus_{i=1}^{\ell+1} e_i\mathcal{U}(G)
= \sum_{i=1}^{\ell+1} \dim_{\mathcal{N}(G)} e_i\mathcal{U}(G)
\ge (\ell+1)/\ell > 1
\]
by \eqref{Esamecor},
and we have a contradiction.  Thus $\mathcal{E}(KG)$ has $n$ primitive
central idempotents $e_1,\dots,e_n$ for some $n\in \mathbb{N}$, $n\le l$.
For each $i$, $1\le i \le n$, choose $0\ne\alpha_i \in
e_i\mathcal{E}(KG)$ such that
$\dim_{\mathcal{N}(G)} \alpha_i \mathcal{U}(G)$ is minimal.

Fix $m \in \{1,2,\dots,n\}$.
Since $\ell \dim_{\mathcal{N}(G)} \alpha \mathcal{U}(G) \in
\mathbb{Z}$ for all $\alpha \in \mathcal{E}(KG)$ by
\eqref{Esamecor}, we may choose $g_1,\dots,g_r \in G$ with
$\dim_{\mathcal{N}(G)} (\sum_{i=1}^r g_i\alpha_m\alpha_m^*g_i^{-1})
\mathcal{U}(G)$ maximal.  Note that if $g_{r+1} \in G$, then
\begin{equation*}
(\sum_{i=1}^{r+1} g_i\alpha_m\alpha_m^*g_i^{-1}) \mathcal{U}(G)
\supseteq \sum_{i=1}^r g_i\alpha_m\mathcal{U}(G)\supseteq (\sum_{i=1}^r
    g_i\alpha_m\alpha_m^*g_i^{-1}) \mathcal{U}(G)
\end{equation*}
by Lemma \ref{Luseful}, hence
\[
\dim_{\mathcal{N}(G)} (\sum_{i=1}^{r+1}
g_i\alpha_m\alpha_m^*g_i^{-1}) \mathcal{U}(G) \ge
\dim_{\mathcal{N}(G)} (\sum_{i=1}^r g_i\alpha_m\alpha_m^*g_i^{-1})
\mathcal{U}(G)
\]
and by maximality of
$\dim_{\mathcal{N}(G)} (\sum_{i=1}^r g_i\alpha_m\alpha_m^*g_i^{-1})
\mathcal{U}(G)$, we see that
\[
\dim_{\mathcal{N}(G)} (\sum_{i=1}^{r+1}
g_i\alpha_m\alpha_m^*g_i^{-1}) \mathcal{U}(G) =
\dim_{\mathcal{N}(G)} (\sum_{i=1}^r g_i\alpha_m\alpha_m^*g_i^{-1})
\mathcal{U}(G).
\]
It follows that
\[
(\sum_{i=1}^{r+1} g_i\alpha_m\alpha_m^*g_i^{-1}) \mathcal{U}(G)
=(\sum_{i=1}^r g_i\alpha_m\alpha_m^*g_i^{-1}) \mathcal{U}(G)
\]
and we deduce from Lemma \ref{Luseful} that
$g \alpha_m \mathcal{U}(G) \subseteq
(\sum_{i=1}^r g_i\alpha_m\alpha_m^*g_i^{-1}) \mathcal{U}(G)$ for all
$g \in G$.  Let $f \in \mathcal{U}(G)$ be the unique projection such
that
\[
f \mathcal{U}(G) =\sum_{i=1}^r g_i\alpha_m\alpha_m^*g_i^{-1}
\mathcal{U}(G).
\]
Then $gf\mathcal{U}(G)=\sum
gg_i\alpha_m\alpha_m^*g_i^{-1}\mathcal{U}(G)\subseteq \sum
gg_i\alpha_m\mathcal{U}(G) \subseteq f\mathcal{U}(G)$ for all $g\in
G$, thus $gf\mathcal{U}(G) = f\mathcal{U}(G)$ and
we deduce that $gfg^{-1}\mathcal{U}(G) = f\mathcal{U}(G)$ for all
$g\in G$.  Also $gfg^{-1}$ is also a projection, thus
$gfg^{-1} = f$ for all $g \in G$ and we conclude that $f$ is a
central projection in $\mathcal{E}(KG)$.  Since
$f \ne 0$, $f\mathcal{U}(G) \subseteq e_m\mathcal{U}(G)$ and $e_m$ is
primitive, we conclude that $f=e_m$ and consequently $\sum_{i=1}^r
g_i\alpha_m\mathcal{U}(G) = e_m\mathcal{U}(G)$.
By omitting some of the terms in this sum if necessary, we may assume
that
\begin{equation} \label{Eextended}
\sum_{1 \le i \le r,\ i\ne s} g_i\alpha_m \mathcal{U}(G) \ne
e_m\mathcal{U}(G)
\end{equation}
for all $s$ such that $1 \le s \le r$.
We make the following observation:
\begin{equation} \label{Eobservation}
\text{If } 0 \ne x \in g_s \alpha_m\mathcal{E}(KG), \text{ then }
x \mathcal{U}(G) = g_s \alpha_m \mathcal{U}(G),
\end{equation}
where $1 \le s \le r$.  This is because $0 \ne x
\mathcal{U}(G) \subseteq g_s\alpha_m \mathcal{U}(G)$ and by
minimality of $\dim_{\mathcal{N}(G)} \alpha_m\mathcal{U}(G)$, we
see that $\dim_{\mathcal{N}(G)} x \mathcal{U}(G) =
\dim_{\mathcal{N}(G)} g_s\alpha_m\mathcal{U}(G)$ and consequently
$x \mathcal{U}(G) = g_s\alpha_m\mathcal{U}(G)$.

We claim that $e_m \mathcal{E}(KG) = \bigoplus_{i=1}^r g_i\alpha_m
\mathcal{E}(KG)$.  Set $\sigma = (\sum_{i=1}^r g_i\alpha_m \alpha_m^*
g_i^{-1})$.  Since $\sigma \mathcal{U}(G) = e_m\mathcal{U}(G)$, we
see that
\begin{equation*}
  (\sigma+ (1-e_m))\mathcal{U}(G) \supseteq \sigma\mathcal{U}(G)+(1-e_m)\mathcal{U}(G)
  = e_m\mathcal{U}(G)+(1-e_m)\mathcal{U}(G) = \mathcal{U}(G).
\end{equation*}
Therefore, $\sigma + 1-e_m$ is invertible in
$\mathcal{U}(G)$ and hence $\sigma + 1-e_m$ is invertible in
$\mathcal{E}(KG)$.  Thus
\begin{equation*}
  e_m\sigma \mathcal{E}(KG)= e_m(\sigma+1-e_m)\mathcal{E}(KG)
  =e_m\mathcal{E}(KG).
\end{equation*}
Moreover, $\sigma\mathcal{E}(KG)\subseteq e_m\mathcal{E}(KG)$ and
therefore $e_m\sigma\mathcal{E}(KG)=\sigma\mathcal{E}(KG)$, hence
\begin{equation*}
e_m\mathcal{E}(KG) = \sigma \mathcal{E}(KG)  =
\sum_{i=1}^r g_i\alpha_m \mathcal{E}(KG).
\end{equation*}
If this sum is not direct, then for some $s$ with $1\le s \le r$, we
have $g_s \alpha_m \mathcal{E}(KG) \cap \sum_{i \ne s} g_i\alpha_m
\mathcal{E}(KG) \ne 0$, and without loss of generality we may assume
that $s=1$.

So let $0 \ne x \in g_1\alpha_m\mathcal{E}(KG) \cap
\sum_{i=2}^r g_i\alpha_m\mathcal{E}(KG)$.  Then $0 \ne x
\mathcal{U}(G) \subseteq g_1\alpha_m \mathcal{U}(G)$ and
\eqref{Eobservation} shows that
$x \mathcal{U}(G) = g_1\alpha_m\mathcal{U}(G)$.  It
follows that $g_1\alpha_m\mathcal{U}(G) \subseteq \sum_{i=2}^r
g_i\alpha_m\mathcal{U}(G)$, consequently
\[
\sum_{i=2}^r g_i\alpha_m\mathcal{U}(G) = e_m\mathcal{U}(G),
\]
which contradicts \eqref{Eextended} and our claim is established.

Now we show that $g_1\alpha_m\mathcal{E}(KG)$ is an irreducible
$\mathcal{E}(KG)$-module.  Suppose $0 \ne x \in g_1\alpha_m
\mathcal{E}(KG)$.  Then $x \mathcal{U}(G)
= g_1 \alpha_m\mathcal{U}(G)$ by \eqref{Eobservation} and using Lemma
\ref{Luseful}, we see as before that $xx^* + \sum_{i=2}^r
g_i\alpha_i\alpha_i^*g_i^{-1} + 1-e_m$ is a unit in $\mathcal{U}(G)$
and hence is also a unit in $\mathcal{E}(KG)$.  This proves that
$x\mathcal{E}(KG) = g_1\alpha_m \mathcal{E}(KG)$ and we deduce that
$\mathcal{E}(KG)$ is a finite direct sum of irreducible
$\mathcal{E}(KG)$-modules.  It follows that $\mathcal{E}(KG)$ is a
semisimple Artinian ring.
\end{proof}

\begin{Prop} \label{Ptorsion}
Let $G$ be a group with $\Delta(G)$ finite and let $K$ be a subfield
of $\mathbb{C}$
with $K=\overline{K}$ which contains all $|{\Delta(G)}|$-th roots of
unity, e.g.~$K=\complexs$ or $K$ is the algebraic
closure of $\rationals$ in $\complexs$.
Then $\mathcal{E}(KG) = \mathcal{D}(KG)$.
\end{Prop}
\begin{proof}
If $e$ is a central idempotent in $\mathcal{U}(G)$, then $e \in
\mathcal{N}(\Delta(G))$, in particular $e \in \mathbb{C}G$, and by our
assumption on $K$ even $e\in KG$. The
result follows.
\end{proof}

The following result is well known, but we include a proof.
\begin{Lem} \label{Ltrace}
Let $G$ be a group, let $e$ be a projection in $\mathcal{N}(G)$, and
let $\alpha \in \mathcal{N}(G)$.  Then $\tr(e\alpha\alpha^*e) \le
\tr(\alpha\alpha^*)$.
\end{Lem}
\begin{proof}
Since $\tr(xy) = \tr(yx)$ for all $x,y\in \mathcal{N}(G)$, we see
that $\tr(e\alpha\alpha^*(1-e)) = \tr((1-e)\alpha\alpha^*e) = 0$.
Therefore $\tr(\alpha\alpha^*) = \tr(e\alpha\alpha^*e) +
\tr((1-e)\alpha\alpha^*(1-e))$.  Since $\tr((1-e)\alpha\alpha^*(1-e))
\ge 0$, the result follows.
\end{proof}

\begin{Lem} \label{Lstrong}
Let $G$ be a group, and let $(\alpha_n)$ be a sequence in
$\mathcal{N}(G)$ converging strongly to $\alpha$.  Suppose that
$\ker\alpha = 0$.  Then $\dim_{\mathcal{N}(G)} (\ker \alpha_n)$
converges to 0.
\end{Lem}
\begin{proof}
By the principle of uniform boundedness, $\| \alpha_n\|$ is bounded.
Also by
multiplying everything by a unitary operator if necessary, we may assume that
$\alpha$ is positive.  Then $\alpha_n - \alpha$ converges strongly to
0 and $(\alpha_n - \alpha)^*$ is bounded, hence $(\alpha_n -
\alpha)^*(\alpha_n - \alpha)$ converges
strongly to 0 and in particular
$\lim_{n\to \infty} \tr ((\alpha_n-\alpha)^*(\alpha_n-\alpha)) =
0$.  Let $e_n \in \mathcal{N}(G)$ denote the projection of $\ell^2(G)$
onto $\ker \alpha_n$.  Then $e_n \alpha_n^* = \alpha_n e_n =0$
and using Lemma \ref{Ltrace}, we obtain
\begin{align*}
\tr((\alpha_n - \alpha)^*(\alpha_n - \alpha)) &\ge
\tr (e_n (\alpha_n - \alpha)^*(\alpha_n - \alpha) e_n)\\
&= \tr(e_n \alpha^*\alpha e_n) \ge 0.
\end{align*}
Thus $\lim_{n\to \infty}\tr(e_n \alpha^*\alpha e_n) = 0$.
Suppose by way of contradiction that
$\lim_{n \to \infty} \dim_{\mathcal{N}(G)} (\ker\alpha_n)
\ne 0$.  Then by
taking a subsequence if necessary, we may assume that
$\dim_{\mathcal{N}(G)} (\ker\alpha_n) > \epsilon$ for
some $\epsilon > 0$, for all $n \in \mathbb{N}$.
By considering the spectral family associated to $\alpha^*\alpha$
\cite[Definition 1.68]{Lueck02}, there is a closed $\alpha^*\alpha$-invariant
$\mathcal{N}(G)$-submodule $X$ of $\ell^2(G)$ and a $\delta > 0$
such that $\dim_{\mathcal{N}(G)}(X) > 1-\epsilon/2$ and
$\alpha^*\alpha > \delta$ on $X$.  Because
$\dim_{\mathcal{N}(G)}(X)>1-\epsilon/2$ and
$\dim_{\mathcal{N}(G)}(\ker\alpha_n)>\epsilon$,
we find that $\dim_{\mathcal{N}(G)} (X\cap
\ker\alpha_n) > \epsilon/2$ (use \cite[Theorem 6.7]{Lueck02}).  Let
$f_n$ denote the projection of $\ell^2(G)$ onto $X \cap \ker\alpha_n$,
so $\tr f_n > \epsilon/2$.
Since $\alpha^*\alpha > \delta$ on $X\cap\ker\alpha_n$,
$f_n\alpha^*\alpha f_n\ge \delta f_n$, and because of positivity of $\tr$
we see that
$\tr(f_n \alpha^*\alpha f_n)\ge \tr(\delta f_n) > \delta\epsilon/2$.
Therefore $\tr (e_n\alpha^*\alpha e_n) > \epsilon\delta/2$ by Lemma
\ref{Ltrace}, which shows that $\tr(e_n\alpha^*\alpha e_n)$ does
not converge to 0, and the result follows.
\end{proof}

\begin{Prop} \label{Pmain}
Let $G$ be a group with $\Delta^+(G) = 1$ and let $K$ be a subfield
of $\mathbb{C}$ that is closed under complex conjugation.
Assume that $\lcm(G) = d
\in \mathbb{N}$ and that $G$ satisfies the strong Atiyah conjecture
over $K$.  Then $\mathcal{D}(KG)$ is a $d
\times d$ matrix ring over a skew field.
\end{Prop}

Let $p$ be a prime, let $q$ be the largest power of $p$ that
divides $d$, and let $H \leq G$ with $|H| = q$ (so $H$ is a ``Sylow"
$p$-subgroup of $G$).  Set $e = \frac{1}{q} \sum_{h\in H} h$, a
projection in $\mathbb{Q}H$.  We shall use the center valued von
Neumann dimension $\dim^u$, as defined in \cite[Definition
9.12]{Lueck02}.  Since $\Delta^+(G) = 1$, we see that
$\dim^u (e\mathcal{U}(G)) = 1/q$ and $\dim^u ((1-e)\mathcal{U}(G)) =
(q-1)/q$.  Therefore by \cite[Theorem 9.13(1)]{Lueck02},
\[
(1-e)\mathcal{U}(G) \cong e\mathcal{U}(G)^{q-1}
\]
and we deduce that there exist orthogonal projections
$e = e_1,e_2,\dots,e_q \in \mathcal{U}(G)$
(so $e_ie_j = 0$ for $i\ne j$) such that $\sum_{i=1}^q
e_i = 1$ and $e_i \mathcal{U}(G) \cong e\mathcal{U}(G)$ for all $i$.
By \cite[Exercise 13.15A, p.~76]{Berberian72}, there exist
similarities (that is self adjoint unitaries)
$u_i \in \mathcal{U}(G)$ with $u_1 = 1$ such that
$e_i = u_ieu_i$.  There is a countable subgroup $F$ of $G$ such that
$u_i \in \mathcal{N}(F)$ for all $i$.  By the Kaplansky density
theorem \cite[Corollary, p.~8]{Arveson76}
for each $i$ ($1 \le i \le q$) there exists a sequence
$u_{ij} \in KF$ such that $u_{ij} \to u_i$ as $j\to \infty$
in the strong operator topology in $\mathcal{N}(F)$
with $u_{1j} = 1$ for all $j$.  Set $v_j =
\sum_{i=1}^q u_{ij}e u_{ij}$.  Then $v_j \to \sum_{i=1}^q e_i = 1$
strongly, hence for $1\le i \le q$,
\[
\lim_{j\to \infty} \dim_{\mathcal{N}(F)}(v_j\mathcal{U}(F)) =
\lim_{j\to \infty} \dim_{\mathcal{N}(F)}(u_{ij}\mathcal{U}(F)) =1
\]
by Lemma \ref{Lstrong}.
Now $\dim_{\mathcal{N}(F)}(x\mathcal{U}(F)) =
\dim_{\mathcal{N}(G)}(x\mathcal{U}(G))$ for all $x \in
\mathcal{U}(F)$, consequently
\begin{equation*}
\lim_{j\to \infty} \dim_{\mathcal{N}(G)}
v_j \mathcal{U}(G) = \lim_{j\to \infty} \dim_{\mathcal{N}(G)}
(u_{ij}\mathcal{U}(G)) =1\qquad \text{for }1 \le i \le q,
\end{equation*}
and since by
assumption $G$ satisfies the strong Atiyah conjecture
over $K$, there exists $n \in \mathbb{N}$ such
that $\dim_{\mathcal{N}(G)} v_j\mathcal{U}(G) =
\dim_{\mathcal{N}(G)} (u_{ij}\mathcal{U}(G)) =1$ for
$1 \le i \le q$ for all $j \ge n$,
in particular $\dim_{\mathcal{N}(G)} (v_n\mathcal{U}(G)) =
\dim_{\mathcal{N}(G)} (u_{in}\mathcal{U}(G)) = 1$ and we
conclude that $v_n$ and $u_{in}$ $(1 \le i \le q)$ are
units in $\mathcal{U}(G)$.  Therefore $v_n$ and $u_{in}$
$(1\le i \le q)$ are units in $\mathcal{D}(KG)$
and we deduce that $\sum_{i=1}^q u_{in}
e\mathcal{D}(KG) = \mathcal{D}(KG)$, because
\begin{equation*}
D(KG) = v_nD(KG) =\sum_{i=1}^q(u_{in}eu_{in})D(KG)\subseteq
\sum_{i=1}^q u_{in}e D(KG)\subseteq D(KG).
\end{equation*}
 Since $\dim_{\mathcal{N}(G)}
e\mathcal{U}(G) = 1/q$, we see that $\bigoplus_{i=1}^q
u_{in}e\mathcal{U}(G) = \mathcal{U}(G)$, a direct sum,
and we deduce that
\begin{equation} \label{Emain}
\bigoplus_{i=1}^qu_{in}e \mathcal{D}(KG) = \mathcal{D}(KG),
\end{equation}
also a direct sum.

Now suppose that $\varepsilon$ is a central idempotent in
$\mathcal{C}(\mathcal{D}(KG))$.  We want to
prove that $\varepsilon =0$ or 1, so assume otherwise.
Now $\varepsilon u_{in}e \mathcal{U}(G)
\cong \varepsilon e \mathcal{U}(G)$ for all $i$, which implies that  $\dim_{\mathcal{N}(G)} (\varepsilon \mathcal{U}(G))
= q \dim_{\mathcal{N}(G)} (\varepsilon e\mathcal{U}(G))$. Moreover, because of
the Atiyah conjecture, $d \dim_{\mathcal{N}(G)} (\varepsilon e
\mathcal{U}(G)) \in
\mathbb{Z}$.   These two observations together imply that
$d\dim_{\mathcal{N}(G)} (\varepsilon\mathcal{U}(G)) \in q\mathbb{Z}$.
Since this is true for all primes $p$, it follows that
$\dim_{\mathcal{N}(G)} \varepsilon\mathcal{U}(G) \in \mathbb{Z}$, so
0 and 1 are the only central idempotents of
$\mathcal{C}(\mathcal{D}(KG))$.

Summing up, we have shown that $\mathcal{C}(\mathcal{D}(KG))$
contains no nontrivial central idempotents.  Using Theorem
\ref{Pextended},
we see that $\mathcal{D}(KG)$ is a semisimple Artinian ring with no
nontrivial central idempotents.  Thus $\mathcal{D}(KG)$ is an
$l\times l$ matrix ring over a division ring for some $l \in
\mathbb{N}$. In particular, $\mathcal{D}(KG)$ is the direct sum of $l$
mutually isomorphic $\mathcal{D}(KG)$-submodules, so if $f$ is a
primitive idempotent in $\mathcal{D}(KG)$, we see that
$\dim_{\mathcal{N}(G)}(f \mathcal{U}(G)) = 1/l$.
Furthermore Lemma \ref{LCramer} (or Lemma \ref{Lsame})
show that $l|d$.  On the other hand \eqref{Emain} shows that $q|l$,
for all primes $p$, so $d|e$ and the result follows.

\begin{proof}[Proof of Theorem \ref{Tmain}]
If $G$ satisfies the strong Atiyah conjecture over $K$, then
$\mathcal{D}(KG)$ is a $d \times d$ matrix ring over
a skew field by Proposition \ref{Pmain}.  Conversely suppose
$\mathcal{D}(KG)$ is a $d \times d$ matrix ring over
a skew field $F$. We need to show that if $M$ is a finitely
presented $KG$-module, then
$\lcm(G)\dim_{\mathcal{N}(G)} M \otimes_{KG} \mathcal{U}(G)
\in \mathbb{Z}$.  However
\[
M \otimes_{KG} \mathcal{U}(G) \cong
M \otimes_{KG} \Mat_d(F) \otimes_{\Mat_d(F)} \otimes \mathcal{U}(G),
\]
consequently $(M \otimes_{KG} \mathcal{U}(G))^d$ is a finitely
generated free $\mathcal{U}(G)$-module and we conclude that
$\lcm(G)\dim_{\mathcal{N}(G)} M \otimes_{KG} \mathcal{U}(G)
\in \mathbb{Z}$ as required.
\end{proof}

\section{Embeddings in $*$-regular rings}\label{sec:embeddings}

There are other closures of group rings $KG$ in $\mathcal{U}(G)$
which may be useful, especially when $\lcm(G) = \infty$.
In general the intersection of regular subrings
of a von Neumann regular ring is not regular
\cite[Example 1.10]{Goodearl91}, however we do have
the following result.

\begin{Prop} \label{Pregular}
Let $G$ be a group and let $\{R_i \mid i\in \mathcal{I}\}$
be a collection of $*$-regular subrings of $\mathcal{U}(G)$.
Then $\bigcap_{i\in \mathcal{I}} R_i$ is also a $*$-regular
subring of $\mathcal{U}(G)$.
\end{Prop}
\begin{proof}
Set $S = \bigcap_{i\in \mathcal{I}} R_i$.  Obviously $S$ is a
$*$-subring of $\mathcal{U}(G)$; we need to show that $S$ is
$*$-regular, that is given $s\in S$, there is a projection $e \in S$
such that $sS = eS$.  We note that $\mathcal{D}(R_i,\mathcal{U}(G))
= R_i$ for all $i$.  Indeed if $x
\in R_i$ and $x$ is invertible in $\mathcal{U}(G)$, then $xR_i = eR_i$
where $e$ is a projection in $R_i$, consequently $x\mathcal{U}(G) =
e\mathcal{U}(G)$ and since $x$ is invertible in $\mathcal{U}(G)$, we
must have $e=1$ and we deduce that $xR_i = R_i$.  Similarly $R_ix =
R_i$ and thus $x$ is invertible in $R_i$, so
$\mathcal{D}(R_i,\mathcal{U}(G)) = R_i$ as asserted.  Since $R_i$ is
$*$-regular, for each $i \in \mathcal{I}$, there is a projection $e_i
\in R_i$ such that $e_iR_i = sR_i$.  We now have $e_i \mathcal{U}(G)$
= $e_j\mathcal{U}(G)$ for all $i,j$ and we deduce that $e_i = e_j$ for
all $i,j \in \mathcal{I}$, so there exists $f \in S$ such
that $f = e_i$ for all $i$.  Since $f\mathcal{U}(G)
= s\mathcal{U}(G)$, we see that $fs = s$, so $s \in fS$ and hence
$sS \subseteq fS$.  Thus the result will be proven if we can show
that $ss^*S \supseteq fS$.  By Lemma \ref{Luseful},
\begin{equation*}
(ss^* + (1-f))\mathcal{U}(G) \supseteq (1-f)\mathcal{U}(G)+s\mathcal{U}(G) =
(1-f)\mathcal{U}(G)+f\mathcal{U}(G)= \mathcal{U}(G)
\end{equation*}
and we see that $ss^* + 1-f$ is a unit in $\mathcal{U}(G)$.
Let $t \in \mathcal{U}(G)$
be the inverse of $ss^* + 1 -f$, so
\begin{equation} \label{Eregular}
(ss^* + 1-f)t = 1.
\end{equation}
Since $\mathcal{D}(R_i,\mathcal{U}(G)) = R_i$ for all $i$, we deduce
that $t \in R_i$ for all $i$ and hence $t \in S$.  Moreover $fs = s$
and $f(1-f) = 0$, so if we multiply \eqref{Eregular} on the left by
$f$, we obtain $ss^*t = f$ and the result is proven.
\end{proof}
Thus if $K$ is a subfield of $\mathbb{C}$ that is closed under
complex conjugation and $G$ is any group, then there is a least
subring of $\mathcal{U}(G)$ containing $KG$ that is $*$-regular.

\bibliographystyle{plain}
\bibliography{atiyahartin}
\end{document}